\documentclass[12pt]{amsart}
\renewcommand{\bar}{\overline}
\usepackage{amsmath,amssymb,amsthm,amscd,euscript}
\renewcommand{\bar}{\overline}
\def \r{\mathbb R}
\def \c{\mathbb C}
\def \q{\mathbb Q}
\def \z{\mathbb Z}

\def \A{\mathcal A}

\def \rat{\hbox{Rat}}

\DeclareMathOperator{\conv}{conv} 

\newtheorem{theorem}{Theorem}[section]
\newtheorem{lemma}[theorem]{Lemma}

\newtheorem{proposition}[theorem]{Proposition}

\theoremstyle{remark}
\newtheorem{remark}[theorem]{Remark}
\theoremstyle{definition}
\newtheorem{definition}[theorem]{Definition}
\newtheorem{example}[theorem]{Example}
\newtheorem{problem}{Problem}
\newtheorem{conjecture}[problem]{Conjecture}

\title[Rational approximation of the maximal commutative subgroups
of $GL(n,\r)$]{Rational approximation of the maximal commutative
subgroups of $GL(n,\r)$}
\author[O.~N.~Karpenkov, A.~M.~Vershik]{Oleg N. Karpenkov$^{1)}$, Anatoly M. Vershik$^{2)}$}
\date{14 October 2009}


\thanks{$^{1)}$Supported by RFBR SS-709.2008.1
 and FWF  No.~S09209.}

\thanks{$^{2)}$Supported by NSh-2460.2008.1, RFBR 08-01-00379, and RFBR 09-01-12175 OFI-M}


 \email[Oleg N. Karpenkov]{karpenk@mccme.ru}
\email[Anatoliy M. Vershik]{vershik@pdmi.ras.ru}

\begin{document}
\input epsf

\begin{abstract}

How to find ``best rational approximations'' of maximal
commutative subgroups of $GL(n,\r)$? In this paper we pose and
make first steps in the study of this problem. It contains both
classical problems of Diophantine and simultaneous approximations
as a particular subcases but in general is much wider. We prove
estimates for $n=2$ for both totaly real and complex cases and
write the algorithm to construct best approximations of a fixed
size. In addition we introduce a relation between best
approximations and sails of cones and interpret the result for
totally real subgroups in geometric terms of sails.
\end{abstract}

\maketitle

\tableofcontents

\section*{Introduction: the problem and its relationships}

We pose and investigate a problem of approximation of maximal
commutative subgroups of $GL(n,\r)$ by rational subgroups, or more
geometrically in other words a problem of approximation of
arbitrary simplicial cones in ${\r}^n$ by rational simplicial
cones. This problem is a natural multidimensional generalization
of a problem on rational approximations of real numbers that is
contained in the case of $n=1$. As a particular example it also
contains a simultaneous approximation problem and closely related
to multidimensional generalizations of continued fractions. The
problem of approximation of real spectrum maximal commutative
subgroups has much in common with the problem of approximations of
nondegenerate simplicial cones. This in particular allows to use
methods dealing with multidimensional continued fractions.

\vspace{2mm}

{
\parindent=0cm
{\bf Maximal commutative subgroups.} We consider a Cartan subgroup
of the group $GL(n,\r)$ or maximal abelian semisimple subgroups of
$GL(n,\r)$. Some times it is convenient to consider such subgroup
as the set of all matrices, commuting with given semisimple
element $A \in GL(n,\r)$, i.e., the centralizer $C_{GL(n,\r)}(A)$.
The centralizer is commutative if and only if $A$ has distinct
eigenvalues. So we work with centralizers of ``generic'' matrices.
For the field of real numbers not all Cartan subgroups are
mutually conjugate: the general Cartan subgroup in $GL(n,\r)$ has
$k$ one-dimensional and $l$ two-dimensional minimal eigenspaces
(where $k{+}2l=n$). We will study mainly the Cartan subgroups with
only one-dimensional minimal eigenspaces, which we call "real
Cartan subgroup", but all the definitions are extended to the
general Cartan subgroups of $GL(n,\r)$ and can be extended to the
case of the Cartan subgroup of $GL(n,\Bbb C)$ or more general
semisimple groups. In that case all elements of the Cartan
subgroup has real eigenvalues. }

\vspace{2mm} We will use term "maximal commutative subgroup" or
shortly MCRF, and denote the space of it as ${\frak C}_n$.

{
\parindent=0cm

\vspace{2mm}
{\bf The space of simplicial cones.}
It is convenient to deal with geometric analog of MCRF-subgroups.
Let us describe a relation of real maximal commutative subgroups
and nondegenerate simplicial cones. }

A {\it nondegenerate simplicial cone} in ${\r}^n$ is a conical
convex hull of a set of $n$ unordered linearly-independent
vectors. Further we omit ``nondegenerate'', since we work only
with nondegenerate cones. Together with any simplicial cone $K$
one may study its symmetric with respect to origin cone $-K$. All
further discussions, constructions, notions, and statements are
invariant with respect to the map $x \mapsto -x$ of ${\r}^n$, and
hence they all deal with both cones $K$ and its symmetric one
$-K$. Therefore, we identify the cones $K$ and $-K$ and define
$Simpl_n$ as {\it a space of pairs of symmetric cones}.

There exists a natural $(2^{n-1})$-folded covering of the space
${\frak C}_n$ of all maximal commutative subgroups by the space
$Simpl_n$:
$$
Simpl_n \rightarrow {\frak C}_n
$$
the cones map to the subalgebras whose eigendirections are the
extremal rays of the cones. So for any element of $Simpl_n$ we
have a maximal commutative subgroups.

Therefore, approximation problems, which we discuss below and
which are local problems, can be studied in terms of the groups as
well as in terms of simplicial cones.

A space $Simpl_n$ of all simplicial cones in ${\r}^n$ can be
defined directly with coordinates of cones generators,
nevertheless it is very important to understand this space as a
{\it homogeneous space of the group $GL(n,\r)$} in the following way.

Consider a group  $GL(n,\r), n>1$ of all linear invertible
transformations in ${\r}^n$ with a fixed basis. Take $D_n$ --- the
subgroup of the diagonal matrices  in the chosen basis which have
positive numbers on the diagonal, i.e. a positive part of the
corresponding Cartan subgroup or connected component of the unity
of that subgroup. The elements of this subgroup leaves invariant
each of the $2^n$ of coordinate cones. The left homogeneous space
$GL(n,\r)/D_n$ can be considered as a space of all connected parts
of the Cartan subgroups of the group $GL(n,\r)$. To get a cone (or
actually a pair of symmetric cones $K$ and $-K$) we should add a
symmetric group of coordinate permutations $S_n$ (Weil group)
which is also contained in the normalizer of $D_n$. Denote by
${\hat D}_n$ the skew-product $S_n\rightthreetimes D_n$ of the
symmetric group and the subgroup of diagonal matrices.

\medskip

{\it A homogeneous space
$$
GL(n,\r)/{\hat D}_n
$$
of left conjugacy classes in $GL(n,\r), n>1$ with respect to the
subgroup ${\hat D}_n$ is naturally identified with the space of
all $($pairs of$)$ nondegenerate simplicial cones $Simpl_n$.}

Indeed, the subgroup of $GL(n,\r)$ preserving the positive
coordinate cone $\r_+^n$ as well as its reflection coincides with
the group ${\hat D}_n$, and $GL(n,\r)$ transitively acts on
$Simpl_n$.

\medskip

Notice that it is sometimes convenient to take the group
$SL(n,\r)$ instead of $GL(n,\r)$ (factoring the last by the
subgroups of positive scalar matrices and taking ${\hat D}_n$ as
the subgroup of positive diagonal matrices with unit determinant
in $ SL(n,\r)$:
$$
Simpl_n =SL(n,\r)/\{ {\hat D}_n \cap SL(n,\r)\}
$$

A homogeneous space $Simpl_n, n>1$  is not compact.  This space
admits a transitive right action of the whole group $GL(n,\r)$ and
it possess an essential absolutely continuous measure $\mu_n$, that is quasihomogeneous with
respect of the action. This measure is called M\"obius measure,
it was studied in~\cite{KarMob}. We are mostly interested in the
actions of $SL(n,\z)$ and $SL(n,\q)$ on the space $Simpl_n$ but
not in the action of the whole group $GL(n,\r), n>1$. These
actions are ergodic.

\begin{definition}
Consider a simplicial cone $C \in Simpl_n$. The boundary of the
convex hall of the integer points in this cone without an origin,
i.e.
$$
\partial \Big(\conv \Big\{ C \cap \z^n\setminus (0,\ldots,
0)\Big\}\Big),
$$
is called {\it the sail of the simplicial cone}.
\end{definition}
The space of the simplicial cones could be identified with the space of the
sails of simplicial cones.

{
\parindent=0cm
{\it Remark.} Note that one can consider the sail for other convex
bodies, for instance of the interiors of conics.
}

For the simplest case of $n=2$ a simplicial cone is a convex angle
between two rays on the plane, and the space $Simpl_2$ of all
cones is a two dimensional torus without a diagonal modulo the
involution: $\{S^1 \times S^1\ \smallsetminus Diag\} /\approx$,
where $Diag$ is the diagonal  in $S^1 \times S^1 $ and $\approx$
is a factorization: $(x,y)\approx (y,x)$. Here the points of the
circles $S^1$ are the oriented lines in $\r^2$ that contains
critical rays of the angles, and quasiinvariant measure is the
Lebesque measure. Actually $Simpl_2$ is a M\"obius strip without a
boundary or equivalently a punctured projective plane. The
geometry of the corresponding cone includes a part of the
classical theory of continuous fraction. The sail for $n=2$ is the
boundary of noncompact convex polygon. The two-dimensional case is
tightly connected with classical continued fractions (see in
Section~2).

\vspace{2mm}

{
\parindent=0cm
{\bf The problem of approximations.}
The described relation between simplicial cones and real spectrum
(i.e.~having real eigenvalues, see further) maximal commutative
subgroups in $GL(n,\r)$ preserving the corresponding cones is a
covering (up to an identification of the cone and its central
symmetrical image). Therefore approximations of such subgroups and
approximations of simplicial cones (we speak about this further)
are the same up to the lifting. Recall that we have fixed a system
of coordinates in $\r^n$, and hence we have a special coordinate
simplicial cone $K_0=\r^n_+$ (a hyperoctant). }

\begin{definition} A {\it rational simplicial cone} $($or respectively
a {\it rational commutative subgroup}$)$ is a cone $($a
subgroup$)$ whose all extremal rays  $($eigen-directions$)$
contains points distinct to the origin with all rational
coordinates, actually this implies the existence of points with
all integer coordinates as well.

A simplicial cone $($maximal commutative subgroup$)$ is called
{\it algebraic} if there exists a matrix $g\in
SL(n,\z)$ with distinct eigenvalues whose eigen-directions generates this cone
$($respectively integer matrix whose centralizer in
$SL(n,\r)$ coincides with this subgroup$)$.
\end{definition}

It is clear that the rational cones form the orbit of the
coordinate cone $K_0$ with respect to the group $SL(n,\q)$.

An example of an algebraic simplicial cone is the conical convex
hull of the two eigenvectors of the Fibonacci matrix:
$$g= \left(
  \begin{array}{cc}
    1  & 1 \\
    1  & 0 \\
  \end{array}
\right)
$$

\medskip
\begin{definition}
 Consider some cone $C \in Simpl_n$
and take nonzero linear forms $L_1, \ldots, L_n$ that annulates
the hyperfaces of the cone. A {\it Markoff-Davenport form} is
$$
\Phi_C(x)= \frac{\prod\limits_{k=1}^{n}\big( L_k(x_1,\ldots,x_n) \big)}
{\Delta(L_1, \ldots, L_n)}
$$
where $\Delta(L_1, \ldots, L_n)$ is the volume of the
parallelepiped spanned by $L_k$ for $k=1,\ldots, n$ in the dual
space.
\end{definition}

This form is defined by a cone uniquely up to a sign. Now having
Markoff-Davenport form $\Phi$ one can define distances between two
cones. For two cones $C_1$ and $C_2$ consider two forms
$$
\Phi_{C_1}(v)+ \Phi_{C_2}(v) \quad \hbox{and} \quad \Phi_{C_1}(v)-
\Phi_{C_2}(v).
$$
Take the maximal absolute values of the coefficients of these
forms separately, the minimal of them would be the distance
between $C_1$ and $C_2$. Further in Subsection~1.1 we define
Markoff-Davenport form in a more general situation.

\medskip

Now we are ready to formulate the main problem of approximations:

{
\parindent=0cm
\emph{\textbf{For a given simplicial cone $($or maximal
commutative subgroup of $SL(n,\r)$$)$ find a rational simplicial
cone $($rational maximal commutative real subgroup$)$ that for a
chosen Markoff-Davenport metric is the closest rational simplicial
cone $($subgroup$)$ in some fixed class of rational cones
$($subgroups$)$}.} }

Such classes of rational cones can chosen to be finite classes
including only cones having fixed ``sizes'' of integer points on
their rays (for more information see below in Section~1).

First of all the approximations problem by
rational simplicial cones (subgroup) must be considered for algebraic cones (subgroups).
The most intriguing things are connected with generalization of the beautiful theory
of Markoff-Lagrange spectra~\cite{Mar} and Markoff-Davenport
$n$-ary forms~\cite{Dav1}.

\vspace{2mm}

{
\parindent=0cm
{\bf Relations with theory of multidimensional continued
fractions.}
The problem on approximation of commutative subgroups or
simplicial cones formulated above and studied in this work is
intimately connected with the theory of multidimensional continued
fractions but does not reduce to that. }

The recent work by V.~I.~Arnold ~\cite{Arn1} and the following works by
him~\cite{Arn2}, E.~I.~Korkina~\cite{Kor2},
G.~Lachaud~\cite{LacBook}, J.-O.~Mussafir~\cite{Mou2},
Karpenkov~\cite{Kar1}, etc., revived the interest to one of classical
generalizations of continued fractions theory, considered for the
first time by F.~Klein in~\cite{Kle1}. From geometrical point of
view the generalization deals with {\it sails}. The
classical theory of ordinary continued fractions i.e. theory of
Gauss transformations in algebro-dynamical terms related to the
case $n=2$ was made by R.~L.~Adler and L.~Flatto in~\cite{Adl}.
M.~L.~Kontsevich, Yu.~M.~Suhov in~\cite{Kon} made an improved
version admitting an extension to multidimensional case.
In the work~\cite{Kon} the authors considered the
following approach to these questions: to study the homogeneous
space $SL(n,\r)/SL(n, \z)$, i.e. the space of lattices
in $SL(n,\r)$, and the action of the Cartan subgroup $D_n$ on it.
For $n=2$ this action is reduced to the action of the group
$\r^1$ and as it is known from~\cite{Adl} it is a special
suspension over the Gauss automorphism that lies in a definition
of continued fractions.

One can suppose that the solution of the approximation problem
reduced to the geometry of the sails in the following sense: in
order to find  the best approximation of the cone (equivalently
maximal commutative subgroup) one must find the appropriate basis
of the vectors which belong to the vertices of the sail of this
cone or adjacent cone. Up to now {\it this is an open question}.
The experiments show that it could be not always the case (see for
instance in Example~\ref{antisail}).

Let us show connections of our problem with this geometry.
First of all the space $Simpl_n$ as we had mentioned can be interpreted as the {\it space of sails
of simplicial cones}. Let us compare our approach to the geometry of sails
with \cite{Kon}.

One can think of dynamical systems as of triples: (a space, a
group action, an invariant or quasiinvariant measure). Then
in~\cite{Kon} the authors study the dynamical system
$$
\{SL(n,\r)/SL(n,\z), \quad D_n,\quad \nu_n\}.
$$
i.e. in our terms it is multidimensional suspension
(time here is a Cartan subgroup) in a given or an arbitrary cone.

Our approach to theory of sails is in some sense dual to the
approach of~\cite{Kon}. We consider another dynamical system,
namely, the action of a discrete (noncommutative) group $SL(n,\z)$
(or $SL(n,\q)$) in the space of sails (or equivalently simplicial
cones):
$$
\{Simpl_n(= SL(n,\r)/{\hat D_n}) , \quad  SL(n,\z),\quad \mu_n\}.
$$
Roughly speaking the ``time'' and the subgroup defining the
homogeneous space has been transposed.

Both approaches have their own advantages and limitations.
However the main aim of the current work is not in studying
of multidimensional sails, their statistics and other properties,
but in their applications to approximations.

\vspace{2mm}
{
\parindent=0cm
{\bf More about geometry of sails.}
 The geometry of sails is very interesting by itself.
 One of the essential subjects here is a statistical analysis of their geometric
characteristics with respect to the measure on the space of the
sails $Simpl_n$. For instance, {\it what is the measure of sails
with given properties: say with given number of faces of some
given combinatorial type} (see~\cite{Kon},~\cite{Avd1},
\cite{Avd2}, \cite{KarZam}, \cite{KarMob}). This would generalize
Gauss-Kuzmin theorem (see in~\cite{Kuz}) and some others for
ordinary continued fractions. The work in this direction has just
started and it is not much known now, first theorems on this
subject can be found in~\cite{KarMob}. }

\vspace{2mm}

Faces of different dimensions of a sail were studied
in~\cite{LacBook}, \cite{Mou2}, \cite{GL}, \cite{Kor1},
\cite{KarPyr}. In algebraic cases all faces are polyhedra. It is
also natural to consider the sails in the adjacent hyperoctants.
The important problem here is to study the condition for a
polygonal surface to be a sail form some cone. This problem was
posed by V.~I.~Arnold and was studied in several papers
(\cite{Arn4}, \cite{Arn2}, \cite{Kar1}, \cite{Kar4D},
\cite{KarPyr}, \cite{KarAlg}, \cite{Kor2}, \cite{Kor3}
\cite{LacBook}, \cite{Mou2}). In~\cite{Tsu} H.~Tsuchihashi showed
the relation between sails of cones and cusp singularities,
introducing a new application to toric geometry. This relation is
studied in detail in~\cite{KarTrig} for the two-dimensional case.

\vspace{2mm}
Actually in the study of $Simpl_n$ the other multidimensional
generalizations of continued fractions can be useful. This in
particular includes the considered before convex-geometric
(\cite{Kle1}, \cite{Arn2}, \cite{Kor2}, \cite{LacBook},
\cite{Kar1}) local minima type (\cite{Min}, \cite{Byk}), Voronoy
(\cite{VAlg}, \cite{Buc}), and algorithmic(\cite{Per}, \cite{Sch})
generalizations of continued fractions.

\vspace{2mm}
{
\parindent=0cm
{\bf Connections with limit shape problems.}
Another link of the approximation problem is with so called limit
shape problems. We want only to emphasize here that the problems
like limit shape problems about Young diagrams or convex lattice
polygons (see~\cite{Ver}) can be considered in the simplicial
cones (instead of traditional posing in the hyperoctant
${\z}_n^+$), and in this case the rational approximation of the
cone becomes an important argument. We hope to consider this in
the appropriate place. }

{
\parindent=0cm
{\bf Description of obtained results.}
Let us  briefly describe the results of this work. Apparently the
problem of approximations of arbitrary commutative subgroups in
$SL(n,\r)$ was never stated in such generality.  By the problem of
approximation we mean the problem of finding of best approximation
of a simplicial cone by rational cones (similar to the classical
problem on best approximations of real numbers by rational
numbers). This problem is very complicated already in the case of
$n=2$. That is also applied even to the algebraic cones. We give
several estimates that suggest an idea that best approximations
are not always related to sails or to sails of adjacent cones (see
also in Example~\ref{antisail}). }

First, we show that the classical case of approximations of real numbers by
rational numbers is really one of particular cases of the proposed new approximation
model. In addition we also indicate that simultaneous approximations are also covered
by our approach.

Further we work in general case of $n=2$. We give upper and lower estimates for the discrepancy
between best approximations and original simplicial cones in the
following important case (Theorem~\ref{Lag}):
{\it let $\alpha_1, \alpha_2\in \r$ both have infinite
continued fractions with bounded elements, consider a simplicial cone
bounded by two lines $y=\alpha_1$ and $y=\alpha_2$, then the growth rate
of the best approximation of size $N$ is bounded by $C_1/N^2$ and $C_2/N^2$ while
$N$ tends to infinity.} Then we translate this statement
to the language of sails and their generalizations
(Theorem~\ref{sails}) and finally show an algorithm to construct best approximations
of a fixed size.

\vspace{2mm}

{
\parindent=0cm
{\it Remark.} In this paper we work in a slightly
extended way including commutative subgroups of $SL(n,\r)$ having
complex conjugate eigenvectors as well. This is the main reason
for our choice to use terminology of commutative subgroups
instead of simplicial cones (that are convenient only for the totally real
case).
}

\vspace{2mm}

We conclude the paper with several examples of approximations in
the three-dimensional case, coming from simultaneous
approximations.

\vspace{2mm}

{
\parindent=0cm
 The paper is organized as follows. In Section 1 we give
basic notions and definitions of maximal subgroup approximation theory.
We introduce sizes and discrepancies for the subgroups and define
the notion of ``best approximations'' in our context. In Section~2 we briefly show
how the classical theory of Diophantine approximations is embedded into
theory of subgroup approximations.
}

Further we make first steps to study a general two-dimensional
case. It is rather complicated since we need to approximate an
object defined by four entries of $2\times 2$ matrices that
vary. Hence this case is comparable with a general case of
simultaneous approximations of vectors in $\r^4$. Nevertheless it
is simpler to find the best approximations in the case of
subgroups, especially in special algebraic case when a certain
periodicity of approximations take place. In Section~3 we write estimates
for the quality of best approximations for both hyperbolic and non-hyperbolic
cases of rays whose continued fractions has bounded elements.
This in particular includes an algebraic case. We also show geometric
origins of the bounds in terms of continued fractions for the hyperbolic algebraic case.

Finally in Section~4
we study in a couple examples the case of simultaneous approximations of vectors in $\r^3$
in the frames of subgroup approximations. We test
two algebraic examples coming from totally real and non-totally real
cases.

\section{Rational approximations of MCRF-groups}
In this section we give general definitions and
formulate basic concepts of maximal commutative subgroups approximations.
We recall a definition of a Markoff-Davenport form in Subsection~1.1. Further
in Subsection~1.2 we define rational subgroups and choose ``size'' for them.
We define the distance function (discrepancy) between two subgroups in
Subsection~1.3.

As we have already mentioned we will continue with terminology of
maximal commutative subgroups. In case when we deal with real
spectra subgroups the statements can be directly translated to
the case of simplicial cones.

\subsection{Regular subgroups and Markoff-Davenport forms}

Consider a real space $\r^n$ and fix some coordinate basis in it.
A real operator is called {\it regular} if all its eigenvalues are
distinct (but not necessary real). A maximal commutative subgroup
of $GL(n,\r)$ is said to be {\it regular}, or {\it MCRS-group} for
short, if it contains regular operators.

We say that a one-dimensional complex space is an {\it eigenspace}
of an MCRF-group if it is an eigenspace of one of its regular
operators. Actually any two regular operators of the same
MCRS-group have the same eigenspaces, therefore each MCRF-group
has exactly $n$ distinct eigenspaces.

Consider an arbitrary MCRS-group $\A$ and denote its eigenspaces by
$l_1,\ldots, l_n$. Denote by $L_i$ a nonzero linear form over
$\c^n$ that attains zero values at all vectors of the complex
lines $l_j$ for $j\ne i$. Let $\Delta(L_1, \ldots, L_n)$ be the
determinant of the matrix having in the $k$-th column the
coefficients of the form $L_k$ for $k=1,\ldots, n$ in the dual
basis.

\begin{definition}
We say that the form
$$
\frac{\prod\limits_{k=1}^{n}\big( L_k(x_1,\ldots,x_n) \big)}
{\Delta(L_1, \ldots, L_n)}
$$
is the {\it Markoff-Davenport form} for the MCRS-group $\A$ and
denote it by $\Phi_{\A}$.
\end{definition}

\begin{example}
Consider an MCRS-group containing a Fibonacci operator
$$
\left(
\begin{array}{cc}
1& 1\\
1& 0\\
\end{array}
\right).
$$
Fibonacci operator has two eigenlines
$$
y=-\theta x \quad \hbox{and}\quad y=\theta^{-1}x,
$$
where $\theta$ is the {\it golden ration} $\frac{1+\sqrt{5}}{2}$.
So the Markoff-Davenport form of Fibonacci operator is
$$
\frac{(y+\theta x)(y-\theta^{-1}x)}{\theta-\theta^{-1}}=
\frac{1}{\sqrt{5}}(-x^2+xy+y^2).
$$

\end{example}

A Markoff-Davenport form is uniquely defined by an MCRS-group up to a
sign, since the linear forms $L_i$ are uniquely defined by the
MCRS-group up to multiplication by a scalar and permutations. By
definition any MCRS-group contains a real operator with distinct
roots, therefore all the coefficients of the Markoff-Davenport
form are real.

\begin{remark}
The minima of the absolute values of such forms on the integer
lattice were studied by A.~Markoff in~\cite{Mar} for
two-dimensional case, and further by H.~Davenport in~\cite{Dav1},
\cite{Dav2}, and~\cite{Dav3} for three-dimensional totally real
case. A few three-dimensional totally real examples were
exhoustively studied by A.~D.~Bryuno, V.~I.~Parusnikov (see for
instance in~\cite{BP}).  The first steps in general
multidimensional case were made in paper~\cite{SL3Z}.
\end{remark}

\subsection{Rational subgroups and their sizes}

We start with the following definition.

\begin{definition}
An MCRS-group $\A$ is called {\it rational} if all its eigenspaces
contain {\it Gaussian} vectors, i. e. vectors whose coordinates
are of type $a+Ib$ for integers $a$ and $b$, where $I^2=-1$.
Denote the set of all rational MCRS-groups of dimension $n$ by $\rat_n$.
\end{definition}

\begin{example}\label{ex2}
The following two operators
$$
\begin{array}{l}
\left(
\begin{array}{rr}
0& -1\\
1&0\\
\end{array}
\right)
\quad \hbox{with eigenvectors $(I, 1)$ and $(-I, 1)$},
\\
\left(
\begin{array}{ll}
1& 1\\
4& 1\\
\end{array}
\right)
\quad \hbox{with eigenvectors $(1, 2)$ and $(1, -2)$}
\end{array}
$$
represents rational MCRS-groups (denote them by $\A_i$ and $\A_{ii}$)
with real and complex conjugate eigen-directions.
\end{example}

For a complex vector $v=(a_1{+}Ib_1,\ldots, a_n{+}Ib_n)$ denote by
$|v|$ the norm
$$
\max\limits_{i=1,\ldots,n}\left(\sqrt{a_i^2+b_i^2}\right).
$$
A Gaussian vector is said to be {\it primitive} if all its
coordinates are relatively prime.

Suppose that a complex one-dimensional space has Gaussian vectors,
then the minimal value of the norm $|*|$ for the Gaussian vectors
is attained at primitive Gaussian vectors.

\begin{definition}
Consider a rational MCRS-group $\A$. Let $l_1,\ldots, l_n$ be the
eigenspaces of $\A$. The {\it size} of $\A$ is a real number
$$
\max\limits_{i=1,\ldots,n}\big\{ |v_i| \big| \hbox{$v_i$ -- is a
primitive Gaussian vector in $l_i$} \big\},
$$
we denote it by $\nu(\A)$.
\end{definition}

The sizes of operators in Example~\ref{ex2} are $1$ and $2$
respectively.

\subsection{Discrepancy functional and approximation model}

We are focused mostly on the following approximation problem:
{\it how to approximate an MCRS-group by rational MCRS-groups $($or even
by a certain subset of rational MCRS-groups$)$}?

Let us first define a natural distance between MCRF-groups. Let
$\A_1$ and $\A_2$ be two MCRS-groups. Consider the following two
symmetric bilinear forms
$$
\Phi_{\A_1}(v)+ \Phi_{\A_2}(v) \quad \hbox{and} \quad
\Phi_{\A_1}(v)- \Phi_{\A_2}(v)
$$
for vectors in $\r^n$. Take the maximal absolute values of the
coefficients of these forms (separately). The minimal of these two maximal values
we consider as a distance between $\A_1$ and $\A_2$,
we call it {\it discrepancy} and denote by $\rho(\A_1,\A_2)$.

\vspace{2mm}

Let us calculate the discrepancy between the MCRS-groups of
Example~\ref{ex2}. We have
$$
\big|\Phi_{\A_i}(v)\pm\Phi_{\A_{ii}}(v)\big|=
\left|I\frac{x^2+y^2}{2} \pm\frac{y^2-4x^2}{4}\right|
$$
therefore $ \rho(\A_i,\A_{ii})=\frac{\sqrt{3}}{2}$.

\begin{definition}
Let $\Omega \subset \rat_n$ for a fixed $n$. The problem of {\it
best approximations} of an MCRS-group $\A$ by MCRS-groups in $\Omega$ is
as follows. {\it For a given positive integer $N$ find a rational
MCRS-group $\A_N$ in $\Omega$ with size not exceeding $N$ such
that}
$$
\rho(\A,\A_N)=\min\big\{\rho(\A,\A') \big|\A'\in \Omega,
\nu(\A')\le N \big\}.
$$
\end{definition}

\remark{There are another important classes of MCRS-groups that
contain matrices of $GL(n,\z)$ and $GL(n,\q)$ respectively. The
MCRS-group is said to be {\it algebraic} if it contains regular
operators of $GL(n,\z)$. It is natural to consider approximations
of MCRS-groups by algebraic MCRS-groups, and approximations of
algebraic MCRS-groups by rational MCRS-groups. }

\section{Diophantine approximations and MCRS-group
approximations}\label{classicsubsection}

A classical problem of approximating real numbers by rational
numbers is a particular case of the problem of best approximations
of MCRS-groups.

For a real $\alpha$ denote by $\A[\alpha]$ an MCRS-group of $GL(2,\r)$
defined by the two spaces $x=0$ and $y=\alpha x$. Consider any two
MCRS-groups $\A[{\alpha_1}]$ and $\A[{\alpha_2}]$ with positive $\alpha_1$ and $\alpha_2$
and calculate a discrepancy between them.
$$
\begin{array}{c}
\displaystyle
\Phi_{\A[{\alpha_1}]}-\Phi_{\A[{\alpha_2}]}=
\frac{x(y-\alpha_1x)}{1}-\frac{x(y-\alpha_2x)}{1}=
(\alpha_2-\alpha_1)x^2
\\
\displaystyle
\Phi_{\A[\alpha_1]}+\Phi_{\A[\alpha_2]}=
\frac{x(y-\alpha_1x)}{1}+\frac{x(y-\alpha_2x)}{1}=
2xy-(\alpha_2+\alpha_1)x^2
\end{array}
$$
Since $\alpha_1>0$ and $\alpha_2>0$ we have
$$
\rho(\A[\alpha_1],\A[\alpha_2])=|\alpha_1-\alpha_2|.
$$

Denote by $\Omega_{[0,1]}^\q$ a subset of all $\A[\alpha]$ for
rational $\alpha$ in the segment $[0,1]$.

For any couple of relatively prime integers $(m,n)$ satisfying $0\le \frac{m}{n}\le 1$ we have
$$
\nu\Big(\A\Big[\frac{m}{n}\Big]\Big)=n.
$$

A classical problem of approximations of real numbers by rational
numbers having bounded denominators in our terminology is as
follows.

\begin{theorem}
Consider a real number $\alpha$, $0\le \alpha \le 1$. Let
$[0,a_1,\ldots]$ $($or $[0,a_1,\ldots, a_k]$$)$ be an ordinary
infinite $($finite$)$ continued fraction for $\alpha$. Then the
set of best approximations consists of MCRS-groups $\A[m/n]$
for $m/n=[0,a_1,\ldots, a_{l-1},a_{l}]$ where $l=1,2,\ldots$ $($In
case of finite continued fraction we additionally have
$\A[m/n]$ for $m/n=[0,a_1,\ldots, a_{k-1},a_{k}{-}1]$$)$.
 \qed
\end{theorem}

\section{General approximations in two-dimensional case}

In this section we prove estimates on the quality of best approximations
for MCRS-groups whose eigen-directions are expressed by
continued fractions with bounded denominators.
We study separately the cases of hyperbolic and non-hyperbolic MCRS-groups.
Especially we study geometric interpretation of the bounds in turms of
geometric continued fractions for the algebraic hyperbolic MCRS-groups.

\subsection{Hyperbolic case}
An MCRS-group is called {\it hyperbolic} if it contains a hyperbolic operator (whose
all eigenvalues are all real and pairwise distinct).

\subsubsection{Lagrange estimates for a special case}

In this subsection we prove an analog of Lagrange theorem on the
approximation rate for an MCRS-groups that has eigenspaces defined by
$y=\alpha_1 x$ and $y=\alpha_2 x$ with bounded elements of the
continued fractions for $\alpha_1$ and $\alpha_2$. In particular
this includes all algebraic MCRS-groups. Here we do not consider the
case when one of the eigenspaces is $x=0$, this case was partially
studied in Section~\ref{classicsubsection}.

\begin{theorem}\label{Lag}
Let $\alpha_1$ and $\alpha_2$ be real numbers having infinite
continued fractions with bounded elements. Consider an MCRS-group
$\A$ with eigenspaces $y=\alpha_1 x$ and $y=\alpha_2 x$.  Then
there exist positive constants $C_1$ and $C_2$ such that for any
positive integer $N$ the best approximation $\A_N$ in $\Omega$
satisfies
$$
\frac{C_1}{N^2}< \rho(\A,\A_N) < \frac{C_2}{N^2}.
$$
\end{theorem}

We will start the proof with the following two lemmas.

Denote by $\A_{\delta_1, \delta_2}$ the MCRS-group defined by the
lines $y=(\alpha_i+\delta_i) x$ for $i=1,2$.
\begin{lemma}\label{lemma1}
Consider a positive real number $\varepsilon_1$ such than
$\varepsilon_1<1/|\alpha_1-\alpha_2|$. Suppose that
$\rho(\A,\A_{\delta_1, \delta_2})<\varepsilon_1$ then
$$
\begin{array}{l}
|\delta_1|<
\frac{(1+|\alpha_1|)(\alpha_1-\alpha_2)^2}{|\alpha_2|(1-\varepsilon_1|\alpha_1-\alpha_2|)}\varepsilon_1
\qquad \hbox{and} \qquad
|\delta_2|<
\frac{(1+|\alpha_2|)(\alpha_1-\alpha_2)^2}{|\alpha_1|(1-\varepsilon_1|\alpha_1-\alpha_2|)}\varepsilon_1.
\end{array}
$$
\end{lemma}

\begin{proof}
Let us remind that the Markoff-Davenport form of
$\A_{\delta_1,\delta_2}$ is
$$
\Phi{\A_{\delta_1, \delta_2}}(x,y)=\frac{\big(y-(\alpha_1 +\delta_1)
x\big)\big(y-(\alpha_2+\delta_2)
x\big)}{(\alpha_2+\delta_2)-(\alpha_1+\delta_1)}.
$$

Consider the absolute values of the coefficients at $y^2$ and at
$xy$ for the difference of Markoff-Davenport forms for the
MCRS-groups $\A$ and $\A_{\delta_1,\delta_2}$. By the conditions of
the lemma these coefficients are less then $\varepsilon_1$:
$$
\left|
\frac{\delta_2-\delta_1}{(\alpha_1-\alpha_2)(\alpha_1-\alpha_2+\delta_1-\delta_2)}\right|<\varepsilon_1
\quad \hbox{and} \quad \left|
\frac{\alpha_1\delta_2-\alpha_2\delta_1}{(\alpha_1-\alpha_2)(\alpha_1-\alpha_2+\delta_1-\delta_2)}\right|<\varepsilon_1.
$$
From the first inequality we have:
$$
|\delta_1-\delta_2|<\frac{(\alpha_1-\alpha_2)^2}{1-\varepsilon_1|\alpha_1-\alpha_2|}\varepsilon_1.
$$
The second inequality implies:
$$
|\delta_1|<\frac{|(\alpha_1-\alpha_2)(\alpha_1-\alpha_2+\delta_1-\delta_2)|\varepsilon_1+
|\alpha_1(\delta_1-\delta_2)|}{|\alpha_2|},
$$
and therefore
$$
|\delta_1|<\frac{|\alpha_1-\alpha_2|(|\alpha_1-\alpha_2|+\frac{(\alpha_1-\alpha_2)^2}
{1-\varepsilon_1|\alpha_1-\alpha_2|}\varepsilon_1)\varepsilon_1+
|\alpha_1|\frac{(\alpha_1-\alpha_2)^2}{1-\varepsilon_1|\alpha_1-\alpha_2|}\varepsilon_1}{|\alpha_2|}=
\frac{(1+|\alpha_1|)(\alpha_1-\alpha_2)^2}{|\alpha_2|(1-\varepsilon_1|\alpha_1-\alpha_2|)}\varepsilon_1.
$$

\vspace{1mm}

The inequality for $\delta_2$ is obtained in the same way.
\end{proof}

\begin{lemma}\label{lemma2}
Let $\varepsilon_2$ be a positive real number. Suppose
$|\delta_1|<\varepsilon_2$ and $|\delta_2|<\varepsilon_2$, then
$$
\rho(\A,\A_{\delta_1,\delta_2})
<\frac{\max\Big(2,2(|\alpha_1|+|\alpha_2|),\alpha_1^2{+}\alpha_2^2+|\alpha_1{-}\alpha_2|\varepsilon_2\Big)}
{(|\alpha_1-\alpha_2|)(|\alpha_1-\alpha_2|+2\varepsilon_2)}
\varepsilon_2.
$$
\end{lemma}

\begin{proof} The statement of lemma follows directly form the
estimate of the coefficients for the difference of
Markoff-Davenport forms for the MCRS-groups $\A$ and
$\A_{\delta_1,\delta_2}$.
\end{proof}

{\it Proof of Theorem~\ref{Lag}.} Let us start with the first
inequality. Let $\alpha_1=[a_0,a_1,\ldots]$, and
$m_i/n_i=[a_0,a_1,\ldots, a_i]$. Without loss of generality we
assume that $N>a_0$. Suppose $k$ is the maximal positive integer
for which $m_k\le N$ and $n_k\le N$. Then we have
$$
\begin{array}{c}
\displaystyle \min\left(\left|\alpha_1-\frac{m}{n}\right|\bigg|
|m|{\le} N, |n|{\le} N \right)\ge
\left|\alpha_1-\frac{m_{k+1}}{n_{k+1}}\right|\ge
\frac{1}{n_{k+1}(n_{k+1}+n_{k+2})}\ge\\
\displaystyle
\frac{1}{(a_{k+1}+1)n_k\big((a_{k+1}+1)n_k+(a_{k+1}+1)(a_{k+2}+1)n_k\big)}\ge
\frac{1}{(a_{k+1}+1)^2(a_{k+2}+2)}\cdot\frac{1}{N^2}.
\end{array}
$$
For the second and the third inequalities we refer to~\cite{Khin}.

The same calculations are valid for $\alpha_2$. Hence we get $C_1$
from Lemma~\ref{lemma1}.

\vspace{1mm}

Now we prove the second inequality.
$$
\left|\alpha_1-\frac{m_k}{n_k}\right|<\frac{1}{n_kn_{k+1}}<\frac{a_{k+1}+1}{n^2_{k+1}}<
\frac{(a_{k+1}+1)}{N^2}\max\big(1,(\alpha_1+1)^2\big).
$$
The first inequality is classical and can be found in~\cite{Khin}.
We take maximum in the last inequality for the case of $m_{k+1}>N$
and $n_{k+1}<N$. From conditions of the theorem the set of $a_i$'s
is bounded. Therefore, there exists a constant $C'_{2,1}$ such
that for any $N$ there exists an approximation of $\alpha_1$ of
quality smaller than $C'_{2,1}/N^2$.

The same holds for $\alpha_2$. Therefore, we can apply
Lemma~\ref{lemma2} in order to obtain the constant $C_2$. \qed

\vspace{2mm}

Let us say a few words about the case of unbounded elements of
continued fractions for $\alpha_i$. Take any positive
$\varepsilon$. If the elements of a continued fraction (say for
$\alpha_1$) are growing fast enough than there exists a sequence
$N_i$ for which the approximations $\A_{N_i}$ are of a quality
$\frac{C}{(N_i)^1+\varepsilon}$. We show this in the following
example.

\begin{example}
Let $M$ be a positive integer. Consider
$\alpha_1=[a_0,a_1,\ldots]$, such that $a_0=1$,
$a_n=(n_{k-1})^{M-1}$. Denote $\frac{m_k}{n_k}=[a_0,\ldots,a_k]$.
Let $\alpha_2=0$. Take $N_k=\frac{n_k+n_{k+1}}{2}$. Then there
exists a positive constant $C$ such that for any integer $i$ we
have
$$
\rho(\A,\A_{N_i})\ge \frac{C}{N_i^{1+1/M}}.
$$
\end{example}

\begin{proof}
For any $i$ we have
$$
n_{i+1}\ge a_in_i=n_i^{M-1}n_i=n_i^{M}.
$$
Therefore, the best approximation with denominator and numerator
less than $N_k$ is not better than
$$
\left|\alpha_1-\frac{m_k}{n_k}\right|\ge
\frac{1}{n_k(n_{k+1}+n_k)}\ge \frac{1}{n_{k+1}^M(n_{k+1}+n_k)}\ge
\frac{2^{1+1/M}}{N_k^{1+1/M}}.
$$
Now we apply Lemma~\ref{lemma1} to complete the proof.
\end{proof}

We suspect the existence of {\it badly approximable} MCRS-group $\A$
and a constant $C$ such that there are only finitely many
solutions $N$ of the following equation
$$
\rho(\A,\A_N)\le \frac{C}{N},
$$
like in the case of simultaneous approximations of vectors in
$\r^3$ (see for instance in~\cite{Laga}).

\subsubsection{Periodic sails and best approximations in algebraic case}

Let us show one relation between classical geometry of numbers  (for example see in~\cite{Arn2}) and best
simultaneous approximations.

First we recall the notion of sails. Consider an arbitrary cone $C$ in $\r^2$
with vertex at the origin and boundary rays $r_1$ and $r_2$.
We also suppose that the angle between $r_1$ and $r_2$ is non-zero and less than $\pi$.
Denote the set of all integer points in the closure of the cone except
the origin by $I_{r_1,r_2}$. The {\it sail} of this
cone is the boundary of the convex hull of $I_{r_1,r_2}$.
It is homeomorphic to a line and
contains rays in case of $r_i$ has an integer point distinct to the origin.

\begin{definition}
Define inductively the {$n$-sail} for the cone $C$.

--- let {\it 1-sail} be the sail of $C$.

--- suppose all $k$-sails for $k<k_0$ are defined then let {\it $k_0$-sail}
be

$$
\partial\Big(\conv\Big(I_{r_1,r_2}\setminus \bigcup\limits_{k=1}^{k_0-1} \hbox{$k$-sail} \Big)\Big),
$$
where $\conv(M)$ denote the convex hull of $M$.
\end{definition}
The $k$-sails have the following interesting property.

\begin{proposition}
Consider a cone $C$. The $k$-sail of $C$ is homothetic to the $1$-sail of $C$ and the coefficient of
homothety is $k$.
\qed
\end{proposition}

Now consider an arbitrary MCRS-group. Let $l_1$ and $l_2$ be the two eigenlines for all the
operators of MCRS-group. The union of all four
$k$-sails for the cones defined by the lines $l_1$ and $l_2$
is a {\it $k$-geometric continued fraction} of the MCRS-group.

\vspace{2mm}

Further we proceed with an algebraic case. So a hyperbolic MCRS-group $\A$ contains an
$GL(2,\z)$-operator with distinct eigenvalues. In this case the mentioned operator acts on
a $k$-geometric continued fraction (for any $k$) as a transitive shift. In addition the values
of the function
$$
\Phi_\A(m,n), \quad \hbox{for $m,n\in \z$,}
$$
are contained in the set $\alpha \z$ where the value $\alpha$ is attained at some
point of the 1-geometric continued fraction. The value $\alpha=\alpha(\A)$ is an essential characteristic of
$\A$, it is sometimes called {\it Markoff minima} of the form $\Phi_\A$.

\begin{lemma}\label{1234}
Let an integer point $(m,n)$ be in the $k$-geometric continued fraction
of $\A$. Then
$$
|\Phi_\A(m,n)|\ge k\alpha.
$$
\end{lemma}

\begin{proof}
We use induction.

The statement clearly holds for $k=1$.

Suppose the statement holds for $k=k_0$ let us prove it
for $k=k_0+1$.
From the step of induction we have the following: for any cone the convex hull of real points
$|\Phi_\A(a,b)|= k_0\alpha$ contains the $k_0$-sail of the cone.
From the other hand all integer points with $|\Phi_\A(m,n)|= k_0\alpha$  (if any)
are on the boundary of this convex hull. Hence all of them are in $k_0$-sail, and thus they are
not contained in $(k_0{+}1)$-sail.
\end{proof}

\begin{theorem}\label{sails}
Let $\A$ be an algebraic MCRS-group.  Then
there exists a positive constants $C$ such that for any
positive integer $N$ the following holds. Let the best approximation $\A_N\in \Omega$
be defined by primitive vectors $v_1$ and $v_2$ contained in $k_1$- and $k_2$-geometric
continued fractions respectively, then $k_1,k_2<C$.
\end{theorem}

\begin{proof}
By Lemma~\ref{1234} it is sufficient to prove that the set of values of $|\Phi_\A(v_i)|$ is bounded.

Let $\A$ has eigenlines $y=\alpha_i x$, $i=1,2$.
Notice that
$$
|\Phi_\A(m,n)|=\left|\frac{(m-\alpha_1 n)(m-\alpha_2 n)}{\alpha_1-\alpha_2}\right|=
\left|\frac{m}{n}-\alpha_1\right|\cdot\left|\frac{m-\alpha_2n}{\alpha_1-\alpha_2} n\right|
$$
Let $v_1=(x_1,y_1)$. By Lemma~\ref{lemma1} (without loss of generality we suppose that
$v_1$ corresponds to $\delta_1$ in the lemma)  the first multiplicative is bounded by
$\tilde C/N^2$ for some constant $\tilde C$ that does not depend on $N$.

Hence,
$$
|\Phi_\A(x_1,y_1)|\le \tilde C \left|\frac{y_1^2}{N^2}\cdot\frac{\frac{x_1}{y_1}-\alpha_2 }{\alpha_1-\alpha_2}\right|\le
\tilde C \left|\frac{\frac{x_1}{y_1}-\alpha_2}{\alpha_1-\alpha_2}\right|
$$
Finally, the last expression is uniformly bounded. The same holds for $v_2$.

Therefore, the set of values of $|\Phi_\A(v_i)|$ is bounded.
\end{proof}

\begin{conjecture}
We conjecture that for almost all $N$ the vectors $v_1$ and $v_2$ defining
$\A_N$ are in $1$-geometric continued fraction.
\end{conjecture}

\subsubsection{Technique of calculation of best approximations in
the hyperbolic case}

In this subsection we show a general technique of calculation of
best approximations for an arbitrary MCRS-group $\A$ with eigenspaces
$y=\alpha_1 x$ and $y=\alpha_2 x$ for distinct real numbers
$\alpha_1$ and $\alpha_2$.

\begin{proposition}\label{lemma3}
Let $m$ and $n$ be two integers. Suppose
$|\alpha_1-\frac{m}{n}|<\varepsilon_3$ $($or
$|\alpha_2-\frac{m}{n}|<\varepsilon_3$  respectively$)$, then the
following holds:
$$
\left|\alpha_1-\frac{m}{n}\right|>
\frac{|\alpha_1-\alpha_2|}{|\alpha_1-\alpha_2|+\varepsilon_3}\frac{|\Phi_\A(m,n)|}{n^2}
\quad \left(
\left|\alpha_2-\frac{m}{n}\right|>\frac{1}{|\alpha_1-\alpha_2|+\varepsilon_3}
\frac{|\Phi_\A(m,n)|}{n^2} \right).
$$
\end{proposition}

\begin{proof}
We have
$$
\begin{array}{l}
\left|\alpha_1-\frac{m}{n}\right|=\frac{1}{n}|m-\alpha_1
n|=\frac{1}{n} \frac{|m-\alpha_1 n|(m-\alpha_2 n)}{m-\alpha_2 n}
=\frac{|\Phi_\A(m,n)|}{n^2}\frac{|\alpha_1-\alpha_2|}{|\alpha_1-\alpha_2+(\frac{m}{n}-\alpha_1)|}>
\frac{|\alpha_1-\alpha_2|}{|\alpha_1-\alpha_2|+\varepsilon_3}\frac{|\Phi_\A(m,n)|}{n^2}.
\end{array}
$$
The same holds for the case of the approximations of $\alpha_2$.
\end{proof}

{\bf Procedure of best approximation calculation}.

{\bf 1).} Find best Diophantine approximations of $\alpha_1$ and
$\alpha_2$ using continued fractions in the square $N\times N$.
Suppose for $\alpha_i$ it is $m_i/n_i$, and the following best
approximation is $m'_i/n'_i$.

{\bf 2).} Consider now the MCRS-group $\bar\A$ with invariant
lines $y=\frac{m_i}{n_i}x$. By Lemma~\ref{lemma2} we get an upper
bound for $\rho(\A,\bar\A)$ (where
$\varepsilon_2=\max(1/(n_1n'_1),1/(n_2n'_2))$).

{\bf 3).} Now having the estimate for discrepancy we use
Lemma~\ref{lemma1} to get estimates $C_1$ and $C_2$ for
$\big|\alpha_1-\frac{p_1}{q_1}\big|$ and
$\big|\alpha_2-\frac{p_2}{q_2} \big|$ for the best approximation
of $\A$ with rays $y=\frac{p_1}{q_1}x$ and $y=\frac{p_2}{q_2}x$.

{\bf 4).} By Proposition~\ref{lemma3} we write an estimate for
$\frac{\Phi_\A(p_i,q_i)}{q_i^2}$ for $i=1,2$.

{\bf 5).} Finally we compare the discrepancies for all MCRS-groups
that satisfies the estimates for $\frac{\Phi_\A(k_i,l_i)}{l_i^2}$
obtained in 4).

\begin{example}
Consider an MCRS-group containing Fibonacci matrix:
$$
\left(
\begin{array}{cc}
0 & 1 \\
1 & 1 \\
\end{array}
\right ).
$$
Denote by $F_n$ the $n$-th Fibonacci number.

Consider any integer $N\ge 100$.

{\bf 1).} Consider a positive integer $k$ such that $F_k\le
N<F_{k+1}$ and choose an approximation $\bar\A$ with eigenspaces
$F_{k-1}y-F_kx=0$ and $F_ky+F_{k-1}x=0$. Then
$$
\left|\alpha_1-\frac{F_k}{F_{k-1}}\right|\le 1/(F_{k-1}F_k), \quad
\left|\alpha_1+\frac{F_{k-1}}{F_k}\right|\le 1/(F_{k}F_{k+1})
$$

{\bf 2).} So, $\varepsilon_2=1/(F_{k-1}F_k)<1/(55\cdot 89)$.
Therefore,
$$
\rho(\A_,\A_{\delta_1,\delta_2})
<\frac{\max\Big(2,2\sqrt{5},3+\sqrt{5}/4895\Big)}
{5+\frac{2\sqrt{5}}{4895}} \frac{1}{F_{k-1}F_k}<\frac{2\sqrt{5}}
{5+\frac{2\sqrt{5}}{4895}}\frac{(89/55)^3}{N^2}<\frac{3.79}{N^2}.
$$

{\bf 3).} Hence, by Lemma~\ref{lemma1} we get
($\varepsilon_1<3.79/100^2$):
$$
\begin{array}{l}
|\delta_1|< \frac{80.35}{N^2}
\qquad \hbox{and} \qquad
|\delta_2|< \frac{18.97}{N^2}.
\end{array}
$$

{\bf 4).} The estimates for $\frac{\Phi_\A(p_1,q_1)}{q_1^2}$ and
$\frac{\Phi_\A(p_2,q_2)}{q_2^2}$ for the corresponding rays of best
approximation are as follows.

$$
\frac{|\Phi_\A(m_1,n_1)|}{n_1^2}<\frac{80.65}{N^2},
\quad
\frac{|\Phi_\A(m_2,n_2)|}{n_2^2}<\frac{18.99}{N^2}.
$$

{\bf 5).} Notice that the number of approximations whose
discrepancies we compare in this step is bounded by some constant
not depending on $N$. We have completed the computations for
$N=10^6$, the answer in this case is the matrix with eigenspaces:
$F_{29}y-F_{30}x=0$ and $F_{30}y+F_{29}x=0$.

We conjecture that for the Fibonacci matrix we always get the best
approximation with eigenspaces $F_{k-1}y-F_kx=0$ and
$F_ky+F_{k-1}x=0$.
\end{example}

We conclude this subsection with an example showing that the
continued fractions do not always give best approximations.

\begin{example}\label{antisail}
Consider an operator $A$ with eigenvectors:
$$
v_1=(1,2) \qquad \hbox{and} \qquad v_2=(2,3),
$$
and the corresponding maximal subgroup $\A$. Then there are four different
best approximations of size 1, they have invariant lines defined by the
following couples of vectors:
$$
\begin{array}{c}
\Big(w_1=(1,0), w_2=(1,1)\Big), \quad
\Big(w_1=(1,0), w_2=(1,-1)\Big),
\\
\Big(w_1=(1,0), w_2=(0,1)\Big), \quad \hbox{and} \quad
\Big(w_1=(0,1), w_2=(1,1)\Big).
\end{array}
$$
(the discrepancy between $\A$ and any of them equals $6$). The
continued fraction (or the union of sails) of $A$ contains only
four integer points
$$
(1,2), \quad (2,3), \quad (-1,-2), \quad \hbox{and} \quad (-2,-3).
$$
Therefore the invariant lines of all four best approximations do not contain vectors of
the sail of $A$.
\end{example}

\begin{remark}
Actually, for a generic MCRS-group the best approximation of any size $N>0$ is unique. In the previous
example we have four best approximations since we are approximating MCRS-group defined by vectors
with integer coefficients.
\end{remark}

\subsection{Non-hyperbolic case}

Now we prove similar statements for the complex case.

\subsubsection{Lagrange estimates for a special case}

In this subsection we prove an analog of Lagrange theorem on the
approximation rate for an MCRS-groups that has complex conjugate
eigenspaces defined by $y=(\alpha+I\beta) x$ and
$y=(\alpha-I\beta) x$ with bounded elements of the continued
fractions for $\alpha$ and $\beta$. In particular this includes
all complex algebraic MCRS-groups.

\begin{theorem}\label{Lag2}
Let $\alpha$ and $\beta$ be real numbers having infinite continued
fractions with bounded elements. Consider an MCRS-group $\A$ with
eigenspaces $y=(\alpha+I\beta) x$ and $y=(\alpha-I\beta) x$. Then
there exist positive constants $C_1$ and $C_2$ such that for any
positive integer $N$ the best approximation $\A_N$ in $\Omega$
satisfies
$$
\frac{C_1}{N^2}< \rho(\A,\A_N) < \frac{C_2}{N^2}.
$$
\end{theorem}

We will start the proof with the following two lemmas.

Denote by $\A_{\delta_1, \delta_2}$ the MCRS-group defined by the
lines $y=\big((\alpha +\delta_1)\pm I(\beta+\delta_2)\big) x$ for
$i=1,2$.
\begin{lemma}\label{lemma2_1}
Consider a positive real number $\varepsilon_1$ such than
$\varepsilon_1<\frac{1}{2(1+|\beta|)}$. Suppose that
$\rho(\A,\A_{\delta_1, \delta_2})<\varepsilon_1$ then
$$
\begin{array}{l}
|\delta_1|< \frac{2|\alpha-\beta|\beta^2}{|\alpha
-\beta|-2\varepsilon_1|\beta|(1+|\beta|)}\varepsilon_1
\qquad \hbox{and} \qquad
|\delta_2|< \frac{2(1+|\beta|+|\alpha-\beta|)\beta^2}{|\alpha
-\beta|-2\varepsilon_1|\beta|(1+|\beta|)}\varepsilon_1.
\end{array}
$$
\end{lemma}

\begin{proof}
Consider the absolute values of the coefficients at $y^2$ and at
$xy$ for the difference of Markoff-Davenport forms for the
MCRS-groups $\A$ and $\A_{\delta_1,\delta_2}$. By the conditions of
the lemma these coefficients are less then $\varepsilon_1$:
$$
\left|
\frac{\delta_2-\delta_1}{2\beta(\beta+\delta_2)}\right|<\varepsilon_1
\quad \hbox{and} \quad \left|
\frac{\alpha\delta_2-\beta\delta_1}{2\beta(\beta+\delta_2)}\right|<\varepsilon_1.
$$
Hence we have
$$
\left|\frac{(\alpha-\beta)\delta_2}{2\beta(\beta+\delta_2)}\right|
\le
+\left|
\frac{\alpha\delta_2-\beta\delta_1}{2\beta(\beta+\delta_2)}\right|
+|\beta|\left|\frac{\delta_2-\delta_1}{2\beta(\beta+\delta_2)}\right|
<(1+|\beta|)\varepsilon_1.
$$
This gives us the estimate for $\delta_2$.

For $\delta_1$ we have
$$
\begin{array}{l}
|\delta_1|<
2|\beta|\left||\beta|+\frac{2(1+|\beta|)\beta^2}{|\alpha
-\beta|-2\varepsilon_1|\beta|(1+|\beta|)}\varepsilon_1\right|
\varepsilon_1+
\frac{2(1+|\beta|)\beta^2}{|\alpha
-\beta|-2\varepsilon_1|\beta|(1+|\beta|)}\varepsilon_1
= \frac{2(1+|\beta|+|\alpha-\beta|)\beta^2}{|\alpha
-\beta|-2\varepsilon_1|\beta|(1+|\beta|)}\varepsilon_1.
\end{array}
$$
The proof is completed.
\end{proof}

\begin{lemma}\label{lemma2_2}
Let $\varepsilon_2$ be a positive real number. Suppose
$|\delta_1|<\varepsilon_2$ and $|\delta_2|<\varepsilon_2$, then
$$
\rho(\A,\A_{\delta_1,\delta_2})
<\frac{\max\Big(2,2(|\alpha|+|\beta|),|\alpha^2{-}\beta^2|+2|\alpha\beta|+2|\beta|\varepsilon_2\Big)}
{|\beta|(|\beta|+\varepsilon_2)} \varepsilon_2.
$$
\end{lemma}

\begin{proof} The statement of lemma follows directly form the
estimate of the coefficients for the difference of
Markoff-Davenport forms for the MCRS-groups $\A$ and
$\A_{\delta_1,\delta_2}$.
\end{proof}

{\it Proof of Theorem~\ref{Lag2}.} The remaining part of the proof
almost completely repeats the end of the proof of
Theorem~\ref{Lag}, so we omit it here. \qed

\subsubsection{Technique of calculation of best approximations in
the hyperbolic case} Here we show a general technique of
calculation of best approximations for an arbitrary MCRS-group $\A$
with eigenspaces $y=(\alpha\pm I\beta)x$ for real number $\alpha$
and positive real $\beta$.

\begin{proposition}\label{lemma2_3}
Let $a$ satisfy $|\alpha+I\beta|<\varepsilon_3$, then the
following holds:
$$
\left|(\alpha+I\beta)-a\right|>
\frac{2\beta|\Phi_\A(1,a)|}{2\beta+\varepsilon_3}.
$$
\end{proposition}

\begin{proof}
We have
$
\begin{array}{l}
\left|(\alpha+I\beta)-a\right|=
\frac{\left|(\alpha+I\beta)-a\right|((\alpha-I\beta)-a)}{(\alpha-I\beta)-a}=
\frac{2\beta|\Phi_\A(1,a)|}{|((\alpha+I\beta)-a)-2I\beta|}>
\frac{2\beta|\Phi_\A(1,a)|}{2\beta+\varepsilon_3}.
\end{array}
$
\end{proof}

{\bf Procedure of best approximation calculation}.

{\bf 1).} Find best Diophantine approximations of $\alpha$ and
$\beta$ using continued fractions in the square $N\times N$.
Suppose for $\alpha$ and $\beta$ it are $m_1/n_1$, and $m_2/n_2$,
and the next best approximation are $m'_1/n'_1$, and $m'_2/n'_2$.

{\bf 2).} Consider the MCRS-group $\bar\A$ with invariant lines
$y=\big(\frac{m_1}{n_1}\pm I\frac{m_2}{n_2}\big)x$. By
Lemma~\ref{lemma2_2} we get an upper bound for $\rho(\A,\bar\A)$
(where $\varepsilon_2=\max(1/(n_1n'_1),1/(n_2n'_2))$).

{\bf 3).} Now having the estimate on discrepancy we use
Lemma~\ref{lemma2_1} to get estimates $C_1$ and $C_2$ for the best
approximation of $\A$: $\big|\alpha-\frac{p_1}{q_1}\big|$ and
$\big|\beta-\frac{p_2}{q_2}\big|$ respectively.

{\bf 4).} By Proposition~\ref{lemma2_3} we write an estimate for
$\big|\Phi_\A\big(1,\frac{p_1}{q_1}+I\frac{p_2}{q_2}\big)\big|$.

{\bf 5).} Finally we compare the discrepancies for all MCRS-groups
that satisfies the estimates obtained in 4).

\section{Simultaneous approximations in $\r^3$ and MCRS-group
approximations}

Theory of simultaneous approximation of a real vector by vectors
with rational coefficients can be considered as a special case of
MCRS-group approximations similarly to the Diophantine case. In this
section we study several examples of simultaneous approximations
in frames of MCRS-group approximations. The first example is an eigen-direction
of a hyperbolic operator (see in Subsection~3.2) and the second is
an eigen-direction of a nonhyperbolic operator (see in Subsection~3.3).

\vspace{2mm}

\subsection{General construction}

Let $[a,b,c]$ be a vector in $\r^3$. Consider the maximal commutative
subgroup $\A[a,b,c]$ defined by three vectors
$$
(a,b,c), \quad (0,1,I), \quad (0,1,-I).
$$

The problem of approximation here is in approximation of the
subgroup $\A[a,b,c]$ by $\A[a',b',c']$ for integer vectors $(a',b',c')$.
For this case we have:
$$
\Phi_{\A[a,b,c]}(x,y,z)=I\left(
-\frac{b^2+c^2}{2a^2}x^3+\frac{b}{a}x^2y+\frac{c}{a}x^2z-\frac{1}{2}
xy^2-\frac{1}{2}xz^2 \right).
$$
Therefore,
$$
\begin{array}{r}
\rho\big(\A[a,b,c],\A[a',b',c']\big)=\min\left(\max\left(\left|\frac{b}{a}-\frac{b'}{a'}\right|,
\left|\frac{c}{a}-\frac{c'}{a'}\right|,
\left|\frac{b^2+c^2}{2a^2}-\frac{{b'}^2+{c'}^2}{2{a'}^2}\right|\right),
\right.\\
\left.
\max \left(\left|\frac{b}{a}+\frac{b'}{a'}\right|,
\left|\frac{c}{a}+\frac{c'}{a'}\right|,
\left|\frac{b^2+c^2}{2a^2}+\frac{{b'}^2+{c'}^2}{2{a'}^2}\right|\right)
\right).
\end{array}
$$

\subsection{A ray of non-hyperbolic operator}

Consider the non-hyperbolic algebraic ope\-rator
$$
B= \left(
\begin{array}{ccc}
0&1&1\\
0&0&1\\
1&0&0\\
\end{array}
\right).
$$
This operator is in some sense the simplest non-hyperbolic
operator we can have (see for more information~\cite{SL3Z}).

Denote the eigenvalues of $E_1$ by $\xi_1$, $\xi_2$, and $\xi_3$
such that $\xi_1$ is real, $\xi_2$ and $\xi_3$ are complex
conjugate. Notice also that
$$
|\xi_1|>|\xi_2|=|\xi_3|.
$$

\vspace{2mm}

We approximate the eigenspace corresponding to $\xi_1$. Let
$v_{\xi_1}$ be the vector in this eigenspace having the first
coordinate equal to 1. Note that
$$
\xi_1\approx 1.3247179573 \quad \hbox{and} \quad
v_{\xi_1}\approx(1, .5698402911, .7548776662).
$$

The set of best approximations $\A_N$ with $N\le 10^6$ contains of
48 elements. These elements are of type $B^{n_i}(1,0,0)$ where
$n_1=4$, and for $2\le i \le 48$ we have $n_i=i+4$. We conjecture
that all the set of best approximations coincide with the set of
points $B^{k}(1,0,0)$ where $k=4$, or $k\ge 6$, the approximation
rate in this case is $CN^{-3/2}$.

\subsection{Two-dimensional golden ratio} Let us consider
an algebraic operator
$$
G= \left(
\begin{array}{ccc}
3&2&1\\
2&2&1\\
1&1&1\\
\end{array}
\right).
$$
This operator is usually called {\it two-dimensional golden
ratio}. It is the simplest hyperbolic operator from many points of
view, his two-dimensional continued fraction in the sense of Klein
was studied in details by E.~I.~Korkina in~\cite{Kor2}
and~\cite{Kor3}.

The group of all integer operators of $GL(3,\z)$ commuting with
$G$ is generated by the following two operators:
$$
E_1= \left(
\begin{array}{ccc}
1&1&1\\
1&1&0\\
1&0&0\\
\end{array}
\right)
\quad \hbox{and} \quad
E_2= \left(
\begin{array}{ccc}
0&1&1\\
1&0&0\\
1&0&-1\\
\end{array}
\right).
$$
Note that $G=E_1^2$ and $E_2=(E_1-Id)^{-1}$, where $Id$ is an
identity operator. Operator $E_1$ is a {\it three-dimensional
Fibonacci operator}.

Denote the eigenvalues of $E_1$ by $\xi_1$, $\xi_2$, and $\xi_3$
in such a way that the following holds:
$$
|\xi_1|>|\xi_2|>|\xi_3|.
$$

\vspace{2mm}

Let us approximate the eigenspace corresponding to $\xi_1$. Denote
by $v_{\xi_1}$ the vector of this eigenspace having the last
coordinate equal to 1. Note that
$$
\xi_1\approx 2.2469796037 \quad \hbox{and} \quad
v_{\xi_1}\approx(2.2469796037, 1.8019377358, 1).
$$

The set of best approximations $\A_N$ with $N\le 10^6$ contains 40
elements. These elements are in the set
$$
\Big\{E_1^mE_2^n(1,0,0)\Big|m,n\in\z\Big\}.
$$

All the points of the sequence can be found from the next table.
In the column $c$ we get $m=m_c$, $n=n_c$ for the approximation
$E_1^{m_c}E_2^{n_c}(1,0,0)$.
\begin{center}

\begin{tabular}{|l||c|c|c|c|c|c|c|c|c|c|c|c|c|c|c|c|c|c|c|c|c|}
\hline
{\bf $i$} & 1 & 2 & 4 & 5 & 6 & 7 & 8 & 9 & 10 & 11
& 12 & 13 & 14 & 15 & 16 & 17 & 18 & 19 & 20 & 21 & 22
\\
\hline \hline
{\bf m}& 1 & 2 & 3 & 3 & 4 & 4 & 5 & 5 & 6 & 6
& 6 & 7 & 7 & 8 & 8 & 9 & 9 & 10 & 10 & 11 & 11
\\
\hline
{\bf n}& 1 & 1 & 2 & 1 & 2 & 1 & 3 & 2 & 3 & 2
& 1 & 3 & 2 & 3 & 2 & 4 & 3 & 4 & 3 & 5 & 4
\\
\hline
\end{tabular}

\vspace{2mm}

\begin{tabular}{|l||c|c|c|c|c|c|c|c|c|c|c|c|c|c|c|c|c|c|c|}
\hline
{\bf $i$} & 23 & 24
& 25 & 26 & 27 & 28 & 29 & 30 & 31 & 32 & 33 & 34
& 35 & 36 & 37 & 38 & 39 & 40 & 41
\\
\hline \hline
{\bf m}& 11 & 12 & 12
& 13 & 13 & 14 & 14 & 15 & 15 & 15 & 16 & 16 & 17
& 17 & 18 & 18 & 19 & 19 & 19
\\
\hline
{\bf n} & 3 & 4 & 3
& 5 & 4 & 5 & 4 & 6 & 5 & 4 & 5 & 4 & 6 & 5 & 6 & 5 & 7
& 6 & 5
\\
\hline
\end{tabular}

\end{center}
In addition to this table we have $\A_3=(3,2,1)$ as best approximation.

We conjecture that all the set of best approximations except
$\A_3$ is contained in the set of all points of type
$E_1^{m}E_2^{n}(1,0,0)$, the approximation rate in this case is $CN^{-3/2}$.


\begin{thebibliography}{99}

\bibitem{Adl}
R.~L.~Adler, L.~Flatto, {\it Cross section maps for geodesic flows.
I. The modular surface.  Ergodic theory and dynamical systems, II}
(College Park, Md., 1979/1980), Progr. Math., v.~21, Birkh\"auser,
Boston, 1982, pp.~103--161.

\bibitem{Arn1}
V.~I.~Arnold, {\it $A$-Graded Algebras and Continued fractions},
Commun. Pure Appl. Math., 142(1989), pp. 993--1000.

\bibitem{Arn4}
V.~I.~Arnold, {\it Higher dimensional continued fractions},
Regular and Chaotic Dynamics, v.~3(3), pp.~10--17, (1998).

\bibitem{Arn2}
V.~I.~Arnold, {\it Continued fractions}, M.: Moscow Center of
Continuous Mathematical Education, 2002.

\bibitem{Avd1}
M.~O.~Avdeeva, V.~A.~Bykovskii, {\it Solution of Arnold's problem
on Gauss-Kuzmin statistics}, Preprint, Vladivostok,
Dal'nauka,~(2002).

\bibitem{Avd2}
M.~O.~Avdeeva, {\it On statistics of incomplete quotients of
finite continued fractions}, Func. an. and appl., v.~38(2004),
n.~2, pp.~1--11.


\bibitem{Byk} V.~A.~Bykovskii, {\it Relative minima of lattices, and vertices of Klein polyhedra},
Funct. Anal. Appl. v.~40(2006),  no.~1, 56--57,

\bibitem{BP}
A.~D.~Bryuno, V.~I.~Parusnikov, {\it Klein polyhedrals for two
cubic Davenport forms}, Mathematical notes, 56(4), (1994),
pp.~9--27.

\bibitem{Buc} J.~A.~Buchmann, {\it A generalization of Voronoi's algorithm I, II}, Journal of Number
Theory, v.~20(1985), pp.~177--209.

\bibitem{Dav1} H.~Davenport, {\it On the product of
three homogeneous linear forms, I}, Proc. London Math. Soc. v.~13(1938),
pp.~139--145.

\bibitem{Dav2} H.~Davenport, {\it Note on the product of
three homogeneous linear forms}, J. London Math. Soc. v.~16(1941),
pp.~98--101.

\bibitem{Dav3} H.~Davenport, {\it On the product of three homogeneous
linear forms. IV}, Math. Proc. Cambridge Philos. Soc., v.~39(1943),
pp~1--21.

\bibitem{GL}
O.~N.~German, E.~L.~Lakshtanov,
{\it On a multidimensional generalization of Lagrange's theorem
for continued fractions} (Russian), Izv. Ross. Akad. Nauk Ser. Mat. v.~(2008),
no.~1, pp.~51--66.


\bibitem{Kar1}
O.~Karpenkov, {\it On tori decompositions associated with
two-dimensional continued fractions of cubic irrationalities},
Func. an. and appl., v.~38(2004), no.~2, pp.~28--37.

\bibitem{KarZam}
O.~Karpenkov, {\it On two-dimensional continued fractions for
integer hyperbolic matrices with small norm}, Russian Math.
Surveys, vol.~59(5), pp.~149--150, 2004.

\bibitem{Kar4D}
O.~Karpenkov, {\it Three examples of three-dimensional continued
fractions in the sense of Klein}, C.~R.~Acad. Sci. Paris, Ser.~B,
v.~343(2006), pp.~5--7.

\bibitem{KarPyr}
O.~Karpenkov, {\it Completely empty pyramids on integer lattices
and two-dimensional faces of multidimensional continued
fractions}, Monatshefte f\"ur Mathematik, v.~152(2007),
pp.~217--249.

\bibitem{KarMob}
O.~Karpenkov,  {\it On invariant M\"obius measure and Gauss-Kuzmin
face distribution}, Proceedings of the Steklov Institute of
Mathematics, v.~258(2007), pp.~74--86.\\
http://arxiv.org/abs/math.NT/0610042.

\bibitem{KarTrig}
O.~Karpenkov, {\it Elementary notions of lattice trigonometry},
Math. Scand., v.~102(2), pp.~161--205, 2008.


\bibitem{KarAlg}
O.~N.~Karpenkov, {\it Constructing multidimensional periodic
continued fractions in the sense of Klein}, Math. Comp. vol.~78,
pp.~1687--1711, 2009.

\bibitem{SL3Z} O.~Karpenkov, {\it Integer conjugacy classes of $SL(3,\z)$
and Hessenberg matrices}, 2007,\\ http://arxiv.org/abs/0711.0830.

\bibitem{Khin}
A.~Ya.~Khinchin, {\it Continued fractions} 4-th ed., ``Nauka'',
Moscow (1978) in Russian; English trans. by Dover Publications,
Inc., Mineola, NY, (1997).
\bibitem{Kle1}
F.~Klein, {\it Ueber einegeometrische Auffassung der gew\"ohnliche
Kettenbruchentwicklung}, Nachr. Ges. Wiss. G\"ottingen Math-Phys.
Kl., v.~3(1891), 357--359.

\bibitem{Kon}
M.~L.~Kontsevich, Yu.~M.~Suhov, {\it Statistics of Klein
Polyhedra and Multidimensional Continued Fractions}, Amer. Math.
Soc. Transl., v.~197(1999), no.~2, pp.~9--27.

\bibitem{Kor1}
E.~I.~Korkina, {\it La p\'eriodicit\'e des fractions continues
multidimensionellez}, C. R. Ac. Sci. Paris, v. 319(1994),
pp.~778--730.

\bibitem{Kor2}
E.~I.~Korkina, {\it Two-dimensional continued fractions. The
simplest examples}, Proceedings of V.~A.~Steklov Math. Ins.,
v.~209(1995), pp.~143--166.

\bibitem{Kor3}
E.~I.~Korkina, {\it The simplest 2-dimensional continued
fraction}, J. Math. Sci., v.~82(1996), no.~5, pp.~3685--3685.

\bibitem{Kuz}
R.~O.~Kuzmin, {\it On a problem of Gauss} Dokl. Akad. Nauk SSSR
Ser A(1928), pp.~375--380.

\bibitem{LacBook}
G.~Lachaud, {\it Voiles et polyh\`edres de Klein},
Act. Sci. Ind., 176~pp, Hermann, 2002.

\bibitem{Laga}
J.~C.~Lagarias,  {\it Best simultaneous Diophantine
approximations. I. Growth rates of best approximation
denominators},  Trans. Amer. Math. Soc. v.~272(1982), no.~2,
pp.~545--554.

\bibitem{Mar} A.~Markoff, {\it Sur les formes quadratiques binaires
ind\'finies}, Math. Ann., v.~15(1879), pp.~381--409.

\bibitem{Min} H.~Minkowski, {\it G\'en\'eralisation de la th\'eorie des
fractions continues}, Ann. Sci. Ecole Norm. Sup., ser.~III, v.~33-4 (1896),
pp.~1057--1070.

\bibitem{Mou2}
J.-O.~Moussafir, {Voiles et Poly\'edres de Klein: Geometrie,
Algorithmes et Statistiques},
docteur en sciences th\'ese, Universit\'e Paris IX - Dauphine, (2000)\\
see also at http://www.ceremade.dauphige.fr/\~{}msfr/

\bibitem{Per} O.~Perron, {\it Grundlagen f\"ur eine Theorie des Jacobischen Kettenbruchalgorithmus},
Math. Ann., v.~64(1907), pp.~1--76.

\bibitem{Sch}
F.~Schweiger,
{\it Multidimensional continued fractions},
Oxford Science Publications, Oxford Univ. Press, viii+234~pp., Oxford, 2000.

\bibitem{Ser}
C.~Series, {\it The modular surface and continued fractions},
J. London Math. Soc. (2), v.~31 (1985), no.~1, 69--80.

\bibitem{Tsu}
H.~Tsuchihashi, {\it Higher dimensional analogues of periodic
continued fractions and cusp singularities}, Tohoku Math. Journ.
v. 35 (1983), pp.~176--193.

\bibitem{Ver}
A.~Vershik. {\it Statistical mechanics of combinatorial partitions, and their limit configurations},
Funct. Anal. Appl. 30, No.2, 90-105 (1996).

\bibitem{VAlg} G.~F.~Voronoy,
{\it On a Generalization of the Algorithm of Continued Fractions},
Izd. Varsh. Univ., Varshava (1896); Collected Works in 3 Volumes
(1952), v.~1, Izd. Akad. Nauk Ukr, SSSR, Kiev (in Russian).
\end{thebibliography}
\end{document}